\documentclass[12pt,letterpaper]{article}
\usepackage[hmargin=1in,vmargin=1in]{geometry}

\usepackage{amsmath}
\usepackage{amsfonts}
\usepackage{amsthm}
\usepackage{amssymb}
\usepackage{array}
\usepackage{caption}
\usepackage{color}
\usepackage{graphics}
\usepackage{graphicx}
\usepackage{hyperref}
\usepackage{verbatim}

\usepackage{tcolorbox}
\tcbuselibrary{theorems}

\newtcbtheorem[]{joke}{}%
{colback=green!5,colframe=green!35!black,fonttitle=\bfseries}{th}

\begin{document}

\title{My Favorite Math Jokes}
\author{Tanya Khovanova}
\maketitle

\begin{abstract}
For many years, I have been collecting math jokes and posting them on my website. I have more than 400 jokes there. In this paper, which is an extended version of my talk at the G4G15, I would like to present 66 of them.
\end{abstract}

\section{Mathematicians and Humor}

Mathematicians are very logical and precise. Here is a joke about that.

\begin{joke*}{}{}
A guy is complaining to his mathematician friend:\\
--- I have a problem. I have difficulty waking up in the morning.\\
--- Logically, counting sheep backwards should help.
\end{joke*}

Did I mention that mathematicians are precise? Once, during a lecture, a professor said the following.

\begin{joke*}{}{}
Assume, for the sake of clarity, that this yellow cube is a blue sphere.
\end{joke*}

Mathematicians are so focused on their abstractions that they ignore the people around them. I met a lot of introverted mathematicians. Here is their favorite joke.

\begin{joke*}{}{}
---What is the difference between an extroverted and an introverted mathematician?\\
---An introverted mathematician will look at his own shoes when he talks to you. An extroverted one will look at your shoes.
\end{joke*}

I was around mathematicians all my life. I even married a mathematician when I was very young. I didn't learn from my mistakes and married a mathematician again, and again.

\begin{joke*}{}{}
---Honey, we are like two parallel lines.\\
---Why do you say that?\\
---The intersection of our life paths was a mistake. 
\end{joke*}

I like humor, but I didn't have enough of it around me. So I started collecting jokes and posting them on my website \cite{TKMH}.

\section{My Own Joke}

Currently, I have more than 400 jokes. Here is a classic one.

\begin{joke*}{}{}
A topologist is someone who doesn't know the difference between a coffee cup and a doughnut.
\end{joke*}

I usually do not invent jokes but rather collect them. Here is one of the exception based on the previous joke.

\begin{joke*}{}{}
A topologist walks into a cafe:\\
--- Can I have a doughnut of coffee, please?
\end{joke*}

Of course, I also told math jokes to my children when they were young; and they too got creative.

\section{Homework Excuses}

My kids heard the standard excuse: my dog ate my homework, and decided to invent more intelligent excuses. They started with a simple variation.

\begin{joke*}{}{}
My biology homework ate my math homework!
\end{joke*}

Then, they made their excuses more mathematical.

\begin{joke*}{}{}
I did part of the homework; the part I have left to do, is 0.999999999...
\end{joke*}

Here is why they didn't do the physics homework.

\begin{joke*}{}{}
I tried to build a black hole in my bedroom when my homework suddenly disappeared.
\end{joke*}

Then, they got interested in computers.

\begin{joke*}{}{}
I accidentally divided by zero and my computer burst into flames.
\end{joke*}

Later, one of them became a programmer.

\begin{joke*}{}{}
My mother redefined my doTheHomework() method with doTheDishes() method.
\end{joke*}

And, as all nerds, they loved ``The Lord of the Rings''.

\begin{joke*}{}{}
My homework was consumed by the power of the One Ring, and it no longer submitted to my will.
\end{joke*}

\section{My Grandson's Joke}

My family consists of nerds. If I asked them, `What's up?', I always get one of two replies, depending on whether we are inside or outside: `the sky' or `the ceiling'. If I ask them to say something, they reply `something'. Now, I tell math jokes to my grandchildren. For example, this famous one.

\begin{joke*}{}{}
A logician rides an elevator. The door opens and someone asks:\\
---Are you going up or down?\\
---Yes.
\end{joke*}

My grandson invented his own silly follow-up for the logicians joke.

\begin{joke*}{}{}
A logician rides an elevator. The door opens and someone asks:\\
---Are you going up or down?\\
---No, replies the logician and walks out.
\end{joke*}

\section{My Teaching}

You won't be surprised to know that the homework I give to my students starts with a math joke. For example, my homework on prime numbers includes the following one.

\begin{joke*}{}{}
Two is the oddest prime.
\end{joke*}

My recent homework on Fibonacci numbers had this joke.

\begin{joke*}{}{}
Fibonacci salad: For today's salad, mix yesterday's leftover salad with that of the day before.
\end{joke*}

Here is one for geometry homework.

\begin{joke*}{}{}
Without geometry, life would be pointless.
\end{joke*}

And another one for algebra.

\begin{joke*}{}{}
---Why was algebra so easy for the Romans?\\
---X was always 10.
\end{joke*}

My topology homework starts with this joke.

\begin{joke*}{}{}
---Why did the chicken cross the M\"{o}bius strip?\\
---To get to the other ... er, um ...
\end{joke*}

I teach a lot of number theory, so I need many jokes about numbers. Here is one of them.

\begin{joke*}{}{}
---Do you know what's odd?\\
---Every other number.
\end{joke*}

Though my students are in middle school, sometimes we dive into calculus.

\begin{joke*}{}{}
---What did the student say about the calculus equation she couldn't solve?\\
---This is derive-ing me crazy!
\end{joke*}

Here is a joke for when I cover statistics.

\begin{joke*}{}{}
---Did you hear about the statistician who drowned crossing a river?\\
---It was three feet deep, on average.
\end{joke*}

I like geometry, and I might give a long lecture about triangles. This is a triangle joke.

\begin{joke*}{}{}
---Why is the obtuse triangle always upset?\\
---It is never right.
\end{joke*}

\section{My Students}

Once, part of the homework I gave was to invent a math joke. Here is what one of my student submitted.

\begin{joke*}{}{}
---Why is Bob scared of the square root of 2?\\
---Because he has irrational fears.
\end{joke*}

Here is another student's creation.

\begin{joke*}{}{}
Everyone envies the circle. It is well-rounded, and highly educated: after all, it has 360 degrees.
\end{joke*}

Yet another student invented a joke on that week's topic: sorting algorithms.

\begin{joke*}{}{}
---Which sorting algorithm is the most relaxing?\\
---The bubble bath sort.
\end{joke*}

\section{Submitted to Me}

Not only do I get new jokes from my family and students, I also receive them from friends and even strangers. This is my friend's daughter talking to her teacher.

\begin{joke*}{}{}
Teacher: What are whole numbers?\\
Student: Like 0, 6, 8, 9.\\
Teacher: And what about 10?\\
Student: It is half-whole, 1 doesn't have a hole.
\end{joke*}

Here is another one from the same source.

\begin{joke*}{}{}
Teacher: Solve the equation: $x + x + x = 9$.\\
Student: $x = 3$, 3, and 3. 
\end{joke*}

This one was received from another friend.

\begin{joke*}{}{}
Quantum entanglement is simple: when you have a pair of socks and you put one of them on your left foot, the other one becomes the ``right sock'', no matter where it is located in the universe. 
\end{joke*}

This joke was submitted by a stranger.

\begin{joke*}{}{}
You have to be odd to be number ONE.
\end{joke*}

\section{Funny Theorems}

Many people laugh at mathematicians. But do mathematician work on funny research? You bet! My coauthor Joel Lewis invented two theorems~\cite{TKTBPT}.

\begin{joke*}{}{}
\textbf{The Big Point Theorem.} Any three lines intersect at a point, provided that the point is big enough.
\end{joke*}

\begin{joke*}{}{}
\textbf{The Thick Line Theorem.} Any three points lie on the same line, provided that the line is thick enough.
\end{joke*}

These ideas can be used for the following famous puzzle. You are given nine dots arranged in three rows and three columns. Connect the dots by drawing four straight, continuous lines that pass through each of the dots without lifting your pencil from the paper. The solution to this puzzle is the most famous example of thinking outside the box. However, there is a humorous solution \cite{TKTIOB} where you can connect the dots using only three lines. Can you figure it out?

\section{Physics}

Mathematics is close to physics. So I have a few physics jokes in my collection. The following joke combines math and physics.

\begin{joke*}{}{}
Einstein-Pythagoras equation: $E = m(a^2 + b^2)$.
\end{joke*}

The next joke is about pure physics.

\begin{joke*}{}{}
Looking for energy?\\
Multiply time by power!
\end{joke*}

Now, a joke about time.

\begin{joke*}{}{}
The barman says, ``We don't serve time travelers in here.''\\
A time traveler walks into a bar.
\end{joke*}

The next joke might be considered to be about time too, but it is about the famous physicist, Werner Heisenberg, and his work.

\begin{joke*}{}{}
Heisenberg gets pulled over on the highway.\\
Cop: ``Do you know how fast you were going, sir?''\\
Heisenberg: ``No, but I know exactly where I am.''
\end{joke*}

After my talk at the G4G15, David Albert sent me a sequel to this joke.

\begin{joke*}{}{}
Heisenberg gets pulled over on the highway.\\
Cop: ``Do you know how fast you were going, sir?''\\
Heisenberg: ``No, but I know exactly where I am.''\\
Cop: ``You were going 85 miles per hour''.\\
Heisenberg: ``Oh great---now I'm lost!''
\end{joke*}

\section{Computer Science}

Mathematics is close to computer science too. Not to mention that one of my children is a computer scientist. So I have a section for computers and programming.

\begin{joke*}{}{}
Humanity invented the decimal system, because people have 10 fingers. And they invented 32-bit computers, because people have 32 teeth.
\end{joke*}

Here is another joke about computers.

\begin{joke*}{}{}
--- Why is your disc drive so noisy?\\
--- It is reading a disc.\\
--- Aloud?
\end{joke*}

And another one about people interacting with computers.

\begin{joke*}{}{}
My computer always beats me in chess. As revenge, I always beat it in a boxing match.
\end{joke*}

The next joke is about computer science's influence on people.

\begin{joke*}{}{}
I saw our system administrator's shopping list. The first line was tomatoes.zip for ketchup.
\end{joke*}

The next joke compares famous modern websites.

\begin{joke*}{}{}
Wikipedia: I know everything.\\
Google: I can find anything.\\
Facebook: I know everyone.\\
Internet: You are nothing without me.\\
Electricity: Shut up, jerks.
\end{joke*}

Here is a recent joke about AI.

\begin{joke*}{}{}
My Roomba has just devoured a piece of cheese I wanted to pick up and eat. The war between humans and robots is already here.
\end{joke*}

\section{Puns}

Because English is my second language, it took me some time to appreciate puns. I started adding them to my collection only recently. Let me start with a self-referencing pun.

\begin{joke*}{}{}
Not all math puns are terrible. Just sum.
\end{joke*}

Speaking about sums\ldots.

\begin{joke*}{}{}
---Why did the two 4's skip lunch?\\
---They already 8!
\end{joke*}

Another pun about numbers.

\begin{joke*}{}{}
---Why shouldn't you argue with a decimal?\\
---Decimals always have a point.
\end{joke*}

Now we move to geometry.

\begin{joke*}{}{}
---Who invented the Round Table?\\
---Sir Cumference.
\end{joke*}

Here is another related one to circumference.

\begin{joke*}{}{}
---How many bakers does it take to bake a pi?\\
---3.14.
\end{joke*}

\section{Teacher Puns}

One category for which I receive a lot of puns is math teachers and schooling. Let's start with multiplication. 

\begin{joke*}{}{}
---Why did the student do multiplication problems on the floor?\\
---The teacher told him not to use tables.
\end{joke*}

Now we go to counting.

\begin{joke*}{}{}
---Why do cheapskates make good math teachers?\\
---Because they make every penny count.
\end{joke*}

Next, we have a geometry teacher joke.

\begin{joke*}{}{}
I saw my math teacher with a piece of graph paper yesterday. I think he must be plotting something.
\end{joke*}

Here is one combining calculus with geometry.

\begin{joke*}{}{}
---Why was math class so long?\\
---The teacher kept going off on a tangent.
\end{joke*}

Here is a joke about behavior of math teachers.

\begin{joke*}{}{}
---What does a hungry math teacher like to eat?\\
---A square meal.
\end{joke*}

I live in Massachusetts, so the next joke is dear to me.

\begin{joke*}{}{}
---What state has the most math teachers?\\
---Math-achusetts.
\end{joke*}

\section{More about Mathematicians}

Now back to mathematicians.

\begin{joke*}{}{}
---Did you hear about the mathematician who's afraid of negative numbers?\\
---He'll stop at nothing to avoid them.
\end{joke*}

Here is a pun about mathematicians.

\begin{joke*}{}{}
---What do mathematicians and the Air Force have in common?\\
---They both use pi-lots.
\end{joke*}

And another one.

\begin{joke*}{}{}
---Where do mathematicians like to go?\\
---Times Square.
\end{joke*}

\section{Kinds of People}

Some jokes go in series. One of my favorite series starts as, ``There are X kinds of people.''

\begin{joke*}{}{}
There are three kinds of people in the world: those who can count and those who can't.
\end{joke*}

The next joke is about binary system.

\begin{joke*}{}{}
There are 10 kinds of people in the world, those who understand binary, and those who don't.
\end{joke*}

Here is a rare fractal joke.

\begin{joke*}{}{}
There are two kinds of people: those who know nothing about fractals and those who think that there are two kinds of people: those who know nothing about fractals and those who think that there are two kinds people\ldots.
\end{joke*}

Data science is adjacent to mathematics. So I include a joke about data here.

\begin{joke*}{}{}
 There are two kinds of people in this world: Those that can extrapolate from incomplete data\ldots.
\end{joke*}

\section{Miscellaneous}

Here are some extra jokes for desert.

\begin{joke*}{}{}
---Did you know that the human brain uses only one third of its capacity?\\
---Hmm, what does the other third do?
\end{joke*}

\begin{joke*}{}{}
---Why was the equal sign so humble?\\
---Because she knew she wasn't greater than or less than anyone else.
\end{joke*}

\begin{joke*}{}{}
---What are ten things you can always count on?\\
---Your fingers.
\end{joke*}

\begin{joke*}{}{}
50\% of marriages end with divorce. The other 50\% end with death.
\end{joke*}

\begin{joke*}{}{}
If a man tries to fail and succeeds, which did he do?
\end{joke*}

\end{document}